\documentclass{article}
\usepackage[a4paper, total={6in, 8in}]{geometry}
\usepackage[utf8]{inputenc}
\usepackage{amssymb}
\usepackage{mathtools}
\usepackage{amsthm}
\usepackage{bbm}
\usepackage{tikz-cd}
\usepackage{dirtytalk}

\newcommand{\bs}[1]{\boldsymbol{#1}}

\title{Smooth relaxation preserving Turing machines}
\author{Adrian K. Xu}
\date{February 2021}
\theoremstyle{definition}
\newtheorem{definition}{Definition}[section]

\theoremstyle{plain}
\newtheorem{theorem}{Theorem}[section]

\theoremstyle{plain}

\theoremstyle{plain}

\theoremstyle{plain}
\newtheorem{lemma}[theorem]{Lemma}

\theoremstyle{remark}
\newtheorem*{remark}{Remark}

\theoremstyle{remark}

\begin{document}

\maketitle

\section{Introduction}

In \cite{dotm}, Clift and Murfet, via ideas from differential linear logic, arrive at a particular means of propagating uncertainty through a Turing machine (TM), interpreted in terms of the state of belief of a \emph{naive Bayesian observer}. The naive Bayesian observer is distinguished from a standard Bayesian observer by a number of independence assumptions--namely, the independence (in distribution) of the tape entries, the state, and the direction of the tape head movement. We henceforth use the phrase \say{naive Bayesian} whenever these assumptions are in effect. We refer to the smooth dynamical system obtained by propagating uncertainty through a TM via naive Bayesian probability as a \emph{smooth relaxation of a TM}.

This leads naturally to a smooth relaxation of the space of programs by considering universal Turing machines (UTM) whose description codes may contain uncertainty, and which are designed in such a way that, when uncertainty is present, their simulation behaviour remains well-defined. Along these lines, in \cite{gops}, the authors introduce the \emph{staged pseudo-UTM}, by which they endow the corresponding subset of Turing machines with a smooth manifold structure. A primary objective of such work is to extend the differential methods which have seen considerable success in modern machine learning to program spaces which more closely approximate the space of computable functions. \\

A number of technical issues regarding the smooth relaxation of TMs remain unclear, a few of which we now survey.

It is well known that the Turing model of computation is robust in the following sense. Given any two variants of the Turing machine model--we may for instance vary the number of tapes or the alphabet size--any partial function computable using one variant is also computable using the other, suitable encodings allowed; typical proofs of results of this kind involve constructing step-wise simulations of arbitrary machines. 

Suppose we generalise the associated partial functions to the smooth case. A formal attempt to do so must resolve the following technicality: what does is mean for the smooth relaxation to halt? Granted a satisfactory answer, we obtain a stronger criterion for the equivalence of two variants of the Turing model--namely, given a machine of one variant, there must exist a machine of the other variant which not only computes the same partial function, but propagates uncertainty in an equivalent manner. A further technicality arises in attempting to formalise this equivalence. It is not always obvious how uncertainty in the tape configuration of one machine can be translated into uncertainty in that of the other machine. 

In this paper, we sidestep the halting technicality by restricting our discussion to step-wise simulations, and by comparing the propagation of uncertainty through two machines in a step-wise manner. For the second technicality, regarding the translation of uncertainty, we for the moment take a pedestrian route and assume uncertainty can be meaningfully translated between two machines using standard probability. These and related issues are elaborated in Sections 2.1 and 2.2. In section 3.1, we show that a multi-tape machine can be simulated on a single tape machine whilst preserving the smooth relaxation, and hence that these two variants are equivalent in the stronger sense just introduced. The equivalence between other variants, in particular, between machines of differing alphabet size, remains to be seen. We shall return to this in the closing remarks.

In section 3.2, we introduce a smooth relaxation preserving pseudo-UTM intended as an alternative to the design presented in \cite{gops}. Conceptually, one would expect a relationship between the way a UTM propagates uncertainty from its simulated input to output, and the particular smooth relaxation of the space of programs it gives rise to. In this vein, we observe at the end of this section that the smooth relaxation preserving pseudo-UTM gives rise to a particularly natural smooth relaxation of the space of programs.

The deeper significance of these technical issues remains unclear. \\

The constructions in this paper can be understood given only a familiarity with Turing machines. However, the formulation and analysis of the smooth relaxation have been expressed in the language of tensor products and direct sums. Hence, at least a passing familiarity with these in the linear algebra context is necessary to appreciate the significance of the constructions.

\section{Preliminaries}

For the present paper, we adopt the following definition of a Turing machine. For our purposes, a given machine need not possess any explicit input-output behaviour, so we omit the initial and halting state from the definition.

\begin{definition}[Turing machine]
A (\emph{single tape}) \emph{Turing machine} (TM) is a triple, $(Q,\Sigma,\delta)$, where $Q$ is a finite set of \emph{states}, $\Sigma$ is a finite \emph{tape alphabet} containing a dedicated \emph{empty symbol} $\square$, and $\delta$ is a \emph{transition function} $Q \times \Sigma \longrightarrow Q \times \Sigma \times \{-1,0,1\}$ sending a \emph{source state} and \emph{read symbol} to a \emph{target state}, \emph{write symbol} and \emph{move direction}.

Denote by $\Sigma^\mathbb{Z}_\square$ the set of functions $\mathbb{Z} \rightarrow \Sigma$ mapping all but finitely many $i \in \mathbb{Z}$ to the blank symbol. Given a transition function (or any function with a Cartesian product codomain) $\delta$, write $\delta_i$ for $\textrm{proj}_i \circ \delta$. Any given Turing machine $(Q,\Sigma,\delta)$ specifies a computation via an associated \emph{step function}, 
\begin{displaymath}
\text{step} : Q \times \Sigma^\mathbb{Z}_\square \longrightarrow Q \times \Sigma^\mathbb{Z}_\square
\end{displaymath}
sending
$$
(q,(y_i)_{i \in \mathbb{Z}}) \longmapsto (q', (y'_{i+d})_{i \in \mathbb{Z}})
$$
where
\begin{align*}
& q' \coloneqq \delta_1(q, y_0) \text{, } \\
& y'_i \coloneqq  
\begin{cases} 
    y_i & i \neq 0 \\
    \delta_2(q, y_0) & i = 0 \\
\end{cases} \text{, and } \\
& d \coloneqq \delta_3(q, y_0).
\end{align*}

\end{definition}

By the \emph{runtime} of a TM on a given input $y \in \Sigma^\mathbb{Z}_\square$ we mean the sequence generated by iteration of the step function on $(q_0,y)$. We refer to an element of this sequence as a \emph{runtime configuration} or just \emph{configuration} and to the second component of such an element as a \emph{tape configuration}.

\subsection{Smooth relaxation}

In the naive Bayesian probabilistic extension, we relax the state to a distribution over the states, and the symbols on the tape to distributions over symbols. We then specify a smooth step function which propagates uncertainty according to the independence assumptions discussed in \cite{dotm}. The language of tensor products provides a natural setting for the constructions to follow, so we formulate our definition in this language and verify (Claim 2.1) that we recover the probability update rules of \cite[Definition~F.1]{gops}.

Given finite set $X$, denote by $\mathbb{R}X$ the free vector space on $X$. Denote by $\Delta X$ the standard $|\Sigma|$-simplex (or rather, its embedding into $\mathbb{R}X$) given by the set $\left\{ \sum_{x \in X} \lambda_x x \in \mathbb{R}X \mid \sum_{x \in X} \lambda_x = 1 \text{ and } \lambda_x \geq 0 \right\}$. Similarly, when $X$ is countably infinite, we denote by $\Delta X$ the set of probability distributions over $X$ with finite support.
Given some function $f: X \rightarrow Y$, we write $\Delta f$ for the unique linear operator $\mathbb{R}X \rightarrow \mathbb{R}Y$ sending each element of the standard basis $x \in X$ to $f(x)$, and we refer to this as the \emph{linear operator induced by} $f$, except where explicitly defined otherwise. 

Note that by $\Delta f_i$ we mean $\Delta (\text{proj}_i \circ f)$. We denote by $\langle \cdot,\cdot \rangle$ the standard inner product.

\begin{definition}[Smooth relaxation of Turing machine]

By the \emph{smooth relaxation of a Turing machine} $(Q,\Sigma,\delta)$ we mean the Turing machine together with a \emph{smooth step function}
$$
\Delta \text{step} : \Delta Q \times (\Delta \Sigma)^\mathbb{Z}_\square \longrightarrow \Delta Q \times (\Delta \Sigma)^\mathbb{Z}_\square
$$
sending
$$
(\bs{q},(\bs{y}_i)_{i \in \mathbb{Z}}) \longmapsto (\bs{q}',(\sum_{d = -1}^1 \langle d, \bs{d} \rangle \bs{y}'_{i+d})_{i \in \mathbb{Z}})
$$
where 
\begin{align*}
& \bs{q}' \coloneqq \Delta \delta_1(\bs{q} \otimes \bs{y}_0) \text{, } \\
& \bs{y}'_i \coloneqq  
\begin{cases} 
    \bs{y}_i & i \neq 0 \\
    \Delta \delta_2(\bs{q} \otimes \bs{y}_0) & i = 0 \\
\end{cases} \text{, and } \\
& \bs{d} \coloneqq \Delta \delta_3(\bs{q} \otimes \bs{y}_0).
\end{align*}

\end{definition}

Here, we view a vector $\bs{q} \in \Delta Q$ as a probability distribution over $Q$ with the probability of a given state $q$ encoded by the scalar projection $\langle q, \bs{q} \rangle$, and similarly for the tape symbols and move direction. As defined above, the operators induced by the transition components, $\Delta \delta_i$, have as their domain $\mathbb{R}(Q \times \Sigma)$, into which the set of probability distributions $\Delta (Q \times \Sigma)$ over $Q \times \Sigma$ embed. Recall there exists an isomorphism $\psi: \mathbb{R}Q \otimes \mathbb{R}\Sigma \xrightarrow{\sim} \mathbb{R}(Q \times \Sigma)$ determined by its operation on the standard bases, $q \otimes \sigma \mapsto (q,\sigma)$. In particular, if $\bs{q} \in \Delta Q$ and $\bs{\sigma} \in \Delta \Sigma$, then $\psi(\bs{q} \otimes \bs{\sigma}) \in \Delta (Q \times \Sigma)$, and we view $\bs{q} \otimes \bs{\sigma}$ as the joint distribution arising when $\bs{q}$ and $\bs{\sigma}$ are treated as the distributions of independent random variables. This framework is employed more generally, and for the most part implicitly, throughout the rest of this paper. \\

Compare the above definition to Definition 2.1. In the conventional situation, the step function would induce an operator $\mathbb{R} (Q \times \Sigma^\mathbb{Z}_\square) \rightarrow \mathbb{R} (Q \times \Sigma^\mathbb{Z}_\square)$, and this would correspond to the propagation of uncertainty according to standard probability. However, here we have pushed in $\mathbb{R}$ (and restricted to the simplices) to obtain the domain and codomain appearing in the definition. This corresponds to the naive Bayesian assumption that the state and each symbol on the tape are conditionally independent at every step.

Once uncertainty has been propagated via $\Delta \delta_i$ from the source state and read symbol to the target state, write symbol and move direction, the resultant tape configuration is computed by a superposition of the tape configuration following the write operation, with one copy for each alternative move direction. The superposition is weighted according to the distribution over the move directions. Note that, according to standard probability, the write symbol distribution $\bs{y}'_{0+d}$ appearing in each term of the superposition should be conditioned on the move direction. The absence of this conditioning corresponds to the further assumption that the move direction and write symbol are conditionally independent. The reader seeking a fuller explanation of this may wish to consult Section 6.2 of \cite{dotm} and compare the equations with the lemma that follows.

Below, we write $\mathbbm{1}()$ for the indicator function, equal to one when the enclosed statement is true and zero otherwise.

The following lemmas follow via direct computation, using linearity of the relevant operators.

\begin{lemma}
In \textnormal{Definition 2.3}, we have the following. Note that the summations are indexed over $Q \times \Sigma$.
\normalfont
\begin{enumerate}
    \item $\langle q_0, \bs{q}' \rangle = \sum_{q, \sigma} \mathbbm{1} (\delta_1(q,\sigma) = q_0) \langle q, \bs{q} \rangle \langle \sigma, \bs{y}_0 \rangle$ for $q_0 \in Q$.
    \item $\langle d_0, \bs{d} \rangle = \sum_{q, \sigma} \mathbbm{1} (\delta_3(q,\sigma) = d_0) \langle q, \bs{q} \rangle \langle \sigma, \bs{y}_0 \rangle$ for $d_0 \in \{-1,0,1\}$.
    \item $\langle \sigma_0, \bs{y}'_0 \rangle = \sum_{q, \sigma} \mathbbm{1} (\delta_2(q,\sigma) = \sigma_0) \langle q, \bs{q} \rangle \langle \sigma, \bs{y}_0 \rangle$ for $\sigma_0 \in \Sigma$.
    \item $\langle \sigma_0, \sum_{d = -1}^1 \langle d, \bs{d} \rangle \bs{y}'_{i+d} \rangle = \sum_{d = -1}^1 \langle d, \bs{d} \rangle [ \mathbbm{1}(i = d) \langle \sigma_0, \bs{y}'_0 \rangle + \mathbbm{1}(i \neq d) \langle \sigma_0, \bs{y}'_{i-d} \rangle ]$ for $\sigma_0 \in \Sigma$.
\end{enumerate}
\end{lemma}

In order to implement the smooth relaxation in our simulator, it is necessary to re-express the tape configuration update rule in the following form.

\begin{lemma}
From \textnormal{Definition 2.3}, we have
$$
\sum_{d = -1}^1 \langle d, \bs{d} \rangle \bs{y}'_{i+d} = \Delta \psi ( \bs{d} \otimes \bs{y}'_{i-1} \otimes \bs{y}'_{i} \otimes \bs{y}'_{i+1} )
$$
where $\psi$ is the function $\{-1,0,1\} \times \Sigma^3 \rightarrow \Sigma$ sending 
$$
(d,\sigma_{-1},\sigma_{0},\sigma_{1}) \mapsto 
\begin{cases} 
\sigma_{-1} & d = -1 \\
\sigma_{0} & d = 0 \\
\sigma_{1} & d = 1 \\
\end{cases}.
$$
It follows that
$$
\sum_{d = -1}^1 \langle d, \bs{d} \rangle \bs{y}'_{i+d} = \Delta \Psi ( \bs{q} \otimes \bs{y}_0 \otimes \bs{y}'_{i-1} \otimes \bs{y}'_{i} \otimes \bs{y}'_{i+1} )
$$
where $\Psi = \psi \circ (\delta_3 \times id)$.
\end{lemma}

One can extend all the above to multi-tape machines as follows. We assume all the tapes share the same alphabet. Recall that an $n$-\emph{tape Turing machine} is given as before by a quadruple, $(Q,\Sigma,\delta,q_0)$, however the transition function is now a map $Q \times \Sigma^n \rightarrow Q \times \Sigma^n \times \{-1,0,1\}^n$, and the step function a map $Q \times (\Sigma^\mathbb{Z}_\square)^n \rightarrow Q \times (\Sigma^\mathbb{Z}_\square)^n$ defined in the obvious manner. In the smooth relaxation, the update rules generalise to:
$$
(\bs{q},(\bs{y}_i^{1})_{i \in \mathbb{Z}},...,(\bs{y}_i^{n})_{i \in \mathbb{Z}}) \longmapsto (\bs{q}',(\sum_{d = -1}^1 \langle d, \bs{d}^{(1)} \rangle \bs{y}'^{(1)}_{i+d})_{i \in \mathbb{Z}}, ..., (\sum_{d = -1}^1 \langle d, \bs{d}^{(n)} \rangle \bs{y}'^{(n)}_{i+d})_{i \in \mathbb{Z}})
$$
where
\begin{align*}
& \bs{q}' \coloneqq \Delta \delta_1(\bs{q} \otimes \bs{y}_0^{(1)} \otimes \cdots \otimes \bs{y}_0^{(n)}) \text{, } \\
& \bs{y}'^{(j)}_i \coloneqq  
\begin{cases} 
    \bs{y}^{(j)}_i & i \neq 0 \\
    \Delta \delta_{1+j}(\bs{q} \otimes \bs{y}_0^{(1)} \otimes \cdots \otimes \bs{y}_0^{(n)}) & i = 0 \\
\end{cases} \text{, for } j=1,...,n \text{, and} \\
& \bs{d}^{(j)} \coloneqq \Delta \delta_{1+n+j}(\bs{q} \otimes \bs{y}_0^{(1)} \otimes \cdots \otimes \bs{y}_0^{(n)}) \text{, for } j=1,...,n.
\end{align*}

From Lemma 2.2 we now have for $j = 1,...,n$

$$
\sum_{d = -1}^1 \langle d, \bs{d}^{(j)} \rangle \bs{y}'^{(j)}_{i+d} = \Delta \Psi^{(j)} ( \bs{q} \otimes \bs{y}_0^{(1)} \otimes \cdots \otimes \bs{y}_0^{(n)} \otimes \bs{y}'^{(j)}_{i-1} \otimes \bs{y}'^{(j)}_{i} \otimes \bs{y}'^{(j)}_{i+1} )
$$
where $\Psi^{(j)} = \psi \circ (\delta_{1+n+j} \times id)$.

\subsection{Smooth relaxation preserving simulations}

Our objective in this section is to make more precise the sense in which our constructions will be smooth relaxation preserving. 

Rogozhin \cite{sutm} defines the notion of simulation as an equivalence between the partial functions computed by each TM. In this paper, we adopt a more narrow viewpoint of a simulation as an equivalence between the dynamical systems specified by the step function of each TM. We will not attempt to formalise this viewpoint; rather, we define a necessary, though not sufficient, criteria for such a simulation, which will provide a minimal context in which the phrase \say{smooth relaxation preserving} has meaning. 

\begin{definition}

Let $M$, $M'$ be Turing machines, with states $Q$, $Q'$, alphabets $\Sigma$, $\Sigma'$, and step functions $step$, $step'$, respectively. Suppose we have a set of \emph{encodings} $Enc \subseteq Q \times \Sigma_\square^\mathbb{Z}$, and a surjective \emph{decoder} $\phi: Enc \rightarrow Q' \times {\Sigma'}_\square^\mathbb{Z}$ such that the following hold.

\begin{enumerate}
\item When run on any input $x \in Q \times \Sigma_\square^\mathbb{Z}$, $M$ passes through $Enc$ infinitely many times.
\item For $x \in Enc$, define $\overline{step}(x) = step^t(x)$, where $t$ is the smallest positive integer such that $step^t(x) \in Enc$. Then the following diagram commutes.

$$
\begin{tikzcd}[row sep = huge, column sep = huge]
Q' \times {\Sigma'}_\square^\mathbb{Z} & Enc \arrow{l}{\phi} \\
Q' \times {\Sigma'}_\square^\mathbb{Z} \arrow{u}{step'} & Enc \arrow{l}{\phi} \arrow{u}{\overline{step}}
\end{tikzcd}
$$

\end{enumerate}

Then we say that the triple $(M,Enc,\phi)$ \emph{generates} $M'$. We shall use the phrase \emph{generating TM} to refer variously to both $M$ and its associated triple. By a \emph{cycle} we mean loosely the sequence of configurations through which $M$ passes when moving from one encoding to the next.

\end{definition}

\begin{remark}
Let $M,M',\phi$ be TMs. If $(M,Enc,\phi)$ generates $M'$ and $(M',Enc',\phi')$ generates $N$, then $(M,\phi^{-1}(Enc'),\phi' \circ \phi)$ generates $N$.
\end{remark}

We expect our notion of a simulation to be a strict subset of the above notion. Indeed, for any TM $M$, there exists a trivial generating TM which stores the initial configuration of $M$ on its tape, alongside a counter, to be incremented on every step; the decoder will then run $M$ on the stored initial configuration for the number of steps recorded by the counter. Such a construction cannot reasonably be characterised as a simulation. Nevertheless, the notion of a generating TM will suffice for our purposes. \\

In the context of Definition 2.3, the smooth configurations $\Delta Q \times (\Delta \Sigma)_\square^\mathbb{Z}$ sit inside the vector space $\mathbb{R} Q \times (\mathbb{R} \Sigma)_\square^\mathbb{Z}$. There exists an embedding, $\iota : \mathbb{R} Q \times (\mathbb{R} \Sigma)_\square^\mathbb{Z} \rightarrow \mathbb{R} (Q \times \Sigma_\square^\mathbb{Z})$ defined as usual by taking the tensor product of the components, such that smooth configurations are sent to their corresponding distribution over classical configurations. That is, the embedding restricts to $\Delta Q \times (\Delta \Sigma)_\square^\mathbb{Z} \rightarrow \Delta (Q \times \Sigma_\square^\mathbb{Z})$. 

We wish to relax the set of classical encodings, $Enc$, to a set of \emph{smooth encodings}, say, $SmoothEnc$. At the very least, we would expect such smooth encodings to be distributions over classical encodings. That is, we expect an inclusion $\iota (SmoothEnc) \subseteq \Delta Enc$. Moreover, we would expect any distribution over classical encodings to be a valid smooth encoding. That is, for any $x$ such that $\iota (x) \in \Delta Enc$, we expect that $x \in SmoothEnc$. Thus, we are obliged to set $SmoothEnc \coloneqq \iota^{-1} (\Delta Enc)$.

A \emph{smooth decoder} must then be a surjective mapping $SmoothEnc \rightarrow \Delta Q' \times (\Delta \Sigma')_\square^\mathbb{Z}$, with appropriate smoothness properties. For now, we shall assume that the uncertainty between the tapes is propagated via standard probability. In this case, the decoder is fully determined by the following composition.
$$
SmoothEnc \xrightarrow{\iota} \Delta Enc \xrightarrow{\Delta \phi} \Delta (Q' \times {\Sigma'}_\square^\mathbb{Z})
$$
Refer to section 2.1 for notational conventions. For this to be a valid decoder, one must verify that its image is indeed $\Delta Q' \times (\Delta \Sigma')_\square^\mathbb{Z}$.

In generalising generating TMs to the smooth setting, we exclude the possibility that the generating TM enters into a configuration which superposes both valid and invalid encodings, hence ensure that the cycles remain cleanly demarcated.

\begin{definition}
Let $(M,Enc,\phi)$ be the generating TM in the context of Definition 2.3. Henceforth, we write $\Delta Enc$ to mean $\iota^{-1} (\Delta Enc)$ and $\Delta \phi$ to mean $\Delta \phi \circ \iota$. We say $(M,Enc,\phi)$ is \emph{well-behaved with respect to the smooth relaxation} if the following hold.
\begin{enumerate}
    \item Given $x \in \Delta Enc$, there exists a $T \in \mathbb{Z}_{> 0}$ such that $(\Delta step)^T (x) \in \Delta Enc$, and $(\Delta step)^t (x) \in \iota^{-1} (\Delta (Q \times \Sigma_\square^\mathbb{Z} \setminus Enc))$ for $0 < t < T$.
    \item $\text{Im } \Delta \phi = \Delta Q' \times (\Delta \Sigma')_\square^\mathbb{Z}$ (relevant embeddings implied). 
\end{enumerate}
We then call $(M,\Delta Enc, \Delta \phi)$ a \emph{smooth generating TM}.
\end{definition}

We are now ready to formalise the meaning of \say{smooth relaxation preserving}.

\begin{definition}
Let $(M,\Delta Enc,\Delta \phi)$ be a smooth generating TM as in Definition 2.4. For $x \in \Delta Enc$, define $\overline{\Delta step}(x) = (\Delta step)^t(x)$, where $t$ is the smallest positive integer such that $(\Delta step)^t(x) \in \Delta Enc$. We say that this smooth generating TM is \emph{smooth relaxation preserving} if the following diagram commutes.

$$
\begin{tikzcd}[row sep = huge, column sep = huge]
\Delta Q' \times (\Delta \Sigma')_\square^\mathbb{Z} & \Delta Enc \arrow{l}{\Delta \phi} \\
\Delta Q' \times (\Delta \Sigma')_\square^\mathbb{Z} \arrow{u}{\Delta step'} & \Delta Enc \arrow{l}{\Delta \phi} \arrow{u}{\overline{\Delta step}}
\end{tikzcd}
$$

\end{definition}

We will henceforth be content to speak semi-formally of \emph{smooth relaxation preserving simulations}, with the understanding that every smooth relaxation preserving simulation is at least a smooth relaxation preserving generating TM. We shall say that a simulation is well-behaved with respect to the smooth relaxation if it is well-behaved as a generating TM. The generalisation to $n$-tape TMs is fairly immediate, and we shall not bother here. \\

Indeed, the staged pseudo UTM given in Appendix F of \cite{gops} is well-behaved with respect to the smooth relaxation, but is not smooth relaxation preserving. For instance, suppose it is initialised with the code for a single state machine which simply writes back to the tape whatever it reads. Suppose the simulation tape alphabet consists of two symbols, $A$ and $B$, and suppose the simulation is initialised with the distribution $0.5\bs{A} + 0.5\bs{B}$ under the tape head. The code contains two tuples, namely, $(A,q,A,q,S)$ and $(B,q,B,q,S)$, where $q$ is the one and only state of the machine. At the beginning of the simulation cycle, the write symbol square on the staging tape reads $X$. The UTM then scans one of the tuples first, say, the one corresponding to $A$. After scanning this tuple, the staged write symbol contains the distribution $0.5\bs{X} + 0.5\bs{A}$. After scanning the second tuple, this distribution becomes $0.5 ( 0.5\bs{X} + 0.5\bs{A} ) + 0.5\bs{B}$, and the final distribution written to the working tape will be $0.375\bs{A} + 0.625\bs{B}$. (The UTM interprets $\bs{X}$ by writing back the read symbol, hence $\bs{X}$ stands in effect for the read distribution $0.5\bs{A} + 0.5\bs{B}$.) This asymmetry does not arise when the uncertainty is propagated directly through the simulated machine. Rather, as one would expect, the distribution remains unchanged.

Moreover, many routine constructions by which multi-tape machines and UTMs are simulated on single-tape machines break down with respect to the smooth relaxation. In general, any construction which relies on certain auxiliary symbols to demarcate segments of the tape and situate the tape head relative to those segments, and which, when passing to the smooth relaxation, experiences ambiguity in the direction of the tape head movement, will not be well-behaved with respect to the smooth relaxation. Inspecting Definition 2.2, we see that any such ambiguity will cause every square along the simulator's tape to be \say{smudged} by the two adjacent squares. That is, uncertainty in a read symbol will, if propagated to uncertainty in the move direction, \say{contaminate} the entire tape. This will typically compromise the simulator's ability to cleanly situate its tape head, and irreversibly distort its working squares.

These considerations are the principle design constraints on the constructions to follow.

\subsection{State partitions and contexts}

In the forthcoming constructions, we view the states as being partitioned into \emph{sections}, each in bijection with a set of \emph{local indices} which we shall call the \emph{context}. The typical context will be a Cartesian product involving the state set and alphabet of the simulated machine. Hence we interpret the simulator's transition function as encoding various components of the simulated machine's transition function.

We shall make this heuristic explicit in our notation, in order to render transparent the behaviour of the construction under its smooth relaxation, and hence its correctness as a smooth relaxation preserving simulation.

Suppose we have a TM with sections $\{Q_i\}_{i \in I}$ and corresponding contexts $\{X_i\}_{i \in I}$. That is, the TM has a set of states $\bigcup_i Q_i$, with each $Q_i$ in bijection with a corresponding set of local indices, $X_i$. Suppose further that there is a family of transitions from $Q_i$ to $Q_j$ for some $i,j \in I$, over some set of read symbols $\Sigma_0 \in \Sigma$. That is, we have $\delta_1 (Q_i \times \Sigma_0) \subseteq Q_j$. Let $f$ be the function $X_i \times \Sigma_0 \rightarrow X_j \times \Sigma \times \{-1,0,1\}$ induced by the restriction $\delta|_{Q_i \times \Sigma_0}$. We refer to $\delta|_{Q_i \times \Sigma_0}$ as the \emph{tract from} $i$ \emph{ to } $j$ \emph{over } $\Sigma_0$, and write the following.
$$
\begin{tikzcd}[column sep = 12em]
[X_i]_i \arrow{r}{\Sigma_0 \rightarrow f_1(x,a), f_2(x,a), f_3(x,a)} & {[X_j]_j}
\end{tikzcd}
$$
Here, the square braces indicate a section with index given by the subscript. The braces enclose the context associated to this section. Hence, states in our designs will never receive explicit names; rather, they will receive a name (an element of $X_i$) local to their section, in terms of which the transition function will be specified. The long arrow denotes the collection of transitions over the states in its source section and an indicated set of read symbols to the left of the smaller arrow (that is, the above mentioned tract). The three expressions to the right of the smaller arrow indicate the target local index in the target section, write symbol, and move direction of an arbitrary transition in the tract in terms of its source local index and read symbol, always denoted by $x$ and $a$ respectively.
\\

We return now to the smooth relaxation. Earlier, we observed that the conditionally independent distributions embed into the set of all distributions over $Q \times \Sigma$ via a map $\Delta Q \times \Delta \Sigma \hookrightarrow \Delta (Q \times \Sigma)$ sending $(\bs{q},\bs{\sigma}) \mapsto \bs{q} \otimes \bs{\sigma}$. Moreover, in Definition 2.2, the linear operator induced by each component of the transition function was treated separately.

By pre-composing with the above embedding the map $\Delta (Q \times \Sigma) \rightarrow \Delta Q \times \Delta \Sigma \times \Delta \{ -1,0,1 \}$ defined component-wise by the linear operator induced by each component of the transition function, we obtain a map $\Delta Q \times \Delta \Sigma \rightarrow \Delta Q \times \Delta \Sigma \times \Delta \{-1,0,1\}$ which we shall call the \emph{smooth transition function}.

Consider the earlier described tract. Suppose the state is distributed over $Q_i$, and the read symbol is distributed over $\Sigma_0$. Then the state in the next time step is distributed over $Q_j$. Hence the restriction of the smooth transition function $\Delta Q_i \times \Delta \Sigma_0 \rightarrow \Delta Q_j \times \Delta \Sigma \times \Delta \{-1,0,1\}$ induces a map $\Delta X_i \times \Delta \Sigma_0 \rightarrow \Delta X_j \times \Delta \Sigma \times \Delta \{-1,0,1\}$ (and vice versa).

In general, in analysing the smooth relaxation of our constructions, we shall view state distributions as distributions over contexts, transformed according the induced map just described. We shall refer to distributions over contexts as \emph{local state distributions}.

\section{Constructions}
\subsection{Multitape on single tape}

\begin{theorem}
There exists a smooth relaxation preserving simulation of any $n$-tape TM on a single tape TM.
\end{theorem}
\begin{proof}
Let $M$ be an $n$-tape TM with states $Q$, alphabet $\Sigma$, and transition function $\delta$. Let $SIM$ be a single tape machine.

The transition function of $SIM$ will be specified in four phases: the \emph{read phase}, \emph{write phase}, \emph{parallel move phase} and \emph{state update phase}. The final section in each phase correspond to the initial section of the next phase, with the final section of the state update phase cycling back to the initial section of the \emph{read phase}. We do not bother to specify unreachable transitions.

The tape alphabet of $SIM$ shall contain $\Sigma$ along with auxiliary symbols $\#_L$, $\#_R$ and $\#_0$.

We now specify an encoding. Let $\bs{y} \in [(\Delta \Sigma)_\square^\mathbb{Z}]^n$ be a tape configuration of $M$. See the end of Section 2.1 for notational conventions. Then an encoding of $\bs{q} \otimes \bs{y}$ in $SIM$ is as follows. $SIM$ shall have a state distribution over section $R1$ (the first section of the read phase) with local distribution $\bs{q}$. Let $R$ be any integer strictly greater than $1$ such that any square on any tape of $M$ with non-zero probability of being non-empty has index less than or equal to $R$. Similarly for $L$. We exclude $-1$, $0$ and $1$ to avoid some edge cases later on. Then $SIM$ shall have a tape configuration containing the sequences given by each row of the following table, interleaved such that $\bs{y}_0^{(1)}$ is at index $0$, $\bs{y}_0^{(2)}$ at index $1$ and so forth.

\begin{center}
\begin{tabular}{ c c c c c c c c c c }
 $\#_L$ & $\bs{y}_L^{(1)}$ & $\cdots$ & $\bs{y}_{-1}^{(1)}$ & $\#_0$ & $\bs{y}_{0}^{(1)}$ & $\bs{y}_{1}^{(1)}$ & $\cdots$ & $\bs{y}_{R}^{(1)}$ & $\#_R$ \\ 
 $\#_L$ & $\bs{y}_L^{(2)}$ & $\cdots$ & $\bs{y}_{-1}^{(2)}$ & $\#_0$ & $\bs{y}_{0}^{(2)}$ & $\bs{y}_{1}^{(2)}$ & $\cdots$ & $\bs{y}_{R}^{(2)}$ & $\#_R$ \\ 
 $\vdots$ &  &  &  &  & \vdots &  &  &  & \vdots \\ 
 $\#_L$ & $\bs{y}_L^{(n)}$ & $\cdots$ & $\bs{y}_{-1}^{(n)}$ & $\#_0$ & $\bs{y}_{0}^{(n)}$ & $\bs{y}_{1}^{(n)}$ & $\cdots$ & $\bs{y}_{R}^{(n)}$ & $\#_R$ \\ 
\end{tabular}
\end{center}

Note that the rightward movement of the tape head amounts to moving down a column and wrapping back up to the top of the next column when the bottom is reached. \\

\underline{Read phase} \\

In the notation of Section 2.3:

$$
\begin{tikzcd}[column sep = 5em]
[Q]_{R1} \arrow{r}{\Sigma \rightarrow (x,a),a,R} & {[Q \times \Sigma]}_{R2} \arrow{r}{\Sigma \rightarrow (x,a),a,R} & \cdots \arrow{r}{\Sigma \rightarrow (x,a),a,R} & {[Q \times \Sigma^{n-1}]}_{Rn} \arrow{r}{\Sigma \rightarrow (x,a),a,S} & {[Q \times \Sigma^n]}_{W1}
\end{tikzcd}
$$

At the beginning of each simulation cycle, $SIM$ will be in the context $Q$, with local state distribution $\bs{q} \in \Delta Q$ mirroring the corresponding state distribution of $M$. The tape head will be over $\bs{y}_0^{(1)}$. Applying the smooth step function, one sees that the tape head moves (unambiguously) right, scanning each encoded read symbol of $MULTI$. The local state distribution transforms as

$$
\bs{q} \mapsto \bs{q} \otimes \bs{y}_0^{(1)} \mapsto \cdots \mapsto \bs{q} \otimes \bs{y}_0^{(1)} \otimes \cdots \otimes \bs{y}_0^{(n)}.
$$ 

In the last tract, the tape head stays put in anticipation of the write phase. \\

\underline{Write phase} \\

Here we write $\delta_{\Sigma,k}$ to denote the component corresponding to the write symbol of the $k$-th tape of $MULTI$.

$$
\begin{tikzcd}[column sep = 8em]
{[Q \times \Sigma^n]}_{W1} \arrow{r}{\Sigma \rightarrow x,\delta_{\Sigma,n}(x),L} &
\cdots \arrow{r}{\Sigma \rightarrow x,\delta_{\Sigma,2}(x),L} & {[Q \times \Sigma^n]}_{Wn} \arrow{r}{\Sigma \rightarrow x,\delta_{\Sigma,1}(x),L} & {[Q \times \Sigma^n]}_{MLB1}
\end{tikzcd}
$$

$SIM$ scans the read symbols in right-to-left, this time replacing them with the appropriate distributions, $\Delta \delta_{\Sigma,k}(\bs{q} \otimes \bs{y}_0^{(1)} \otimes \cdots \otimes \bs{y}_0^{(n)})$. (Refer to the end of Section 2.1 for the smooth step function of multi-tape machines.) \\

\underline{Parallel move phase} \\

This is the phase most sensitive to our requirement that the simulation be smooth relaxation preserving.

In order to compute the superposition at an encoded square of $MULTI$, say at index $i$ of tape $k$, resulting from the possibly ambiguous tape head movement, $SIM$ must have loaded into its local state distribution the symbol distributions at indices $i-1$, $i$ and $i+1$. Since the original symbol distribution on one side is erased by the update operation, $SIM$ must \say{remember} the original symbol distribution of the most recently updated square, then discard it once all superpositions involving it have been computed and written to the tape. An alternative solution would be to introduce staging squares between the working squares to store copies of the original distributions.

We only give a partial construction, which computes the superpositions in the bottom row of the table on page 10. Extending the construction to the remaining rows is routine, and can be achieved by copy-pasting the partial construction with appropriate intermediate transitions and minor adjustments. Alternatively, one would loop a suitably modified partial construction over the rows and introduce new auxiliary symbols to cue the loop exit. As usual, there will be a trade-off between the size of the alphabet and the number of states used. As we have no intention of being economical in this regard, we proceed without further comment.

We divide the partial construction into four sub-phases: \emph{left border shift}, \emph{left edge case}, \emph{main loop}, \emph{right edge case} and \emph{right border shift}. These will be indicated by prefixes $MLB$, $MLE$, $MML$, $MRE$ and $MRB$ in the section indices. \\

\emph{Left border shift} \\

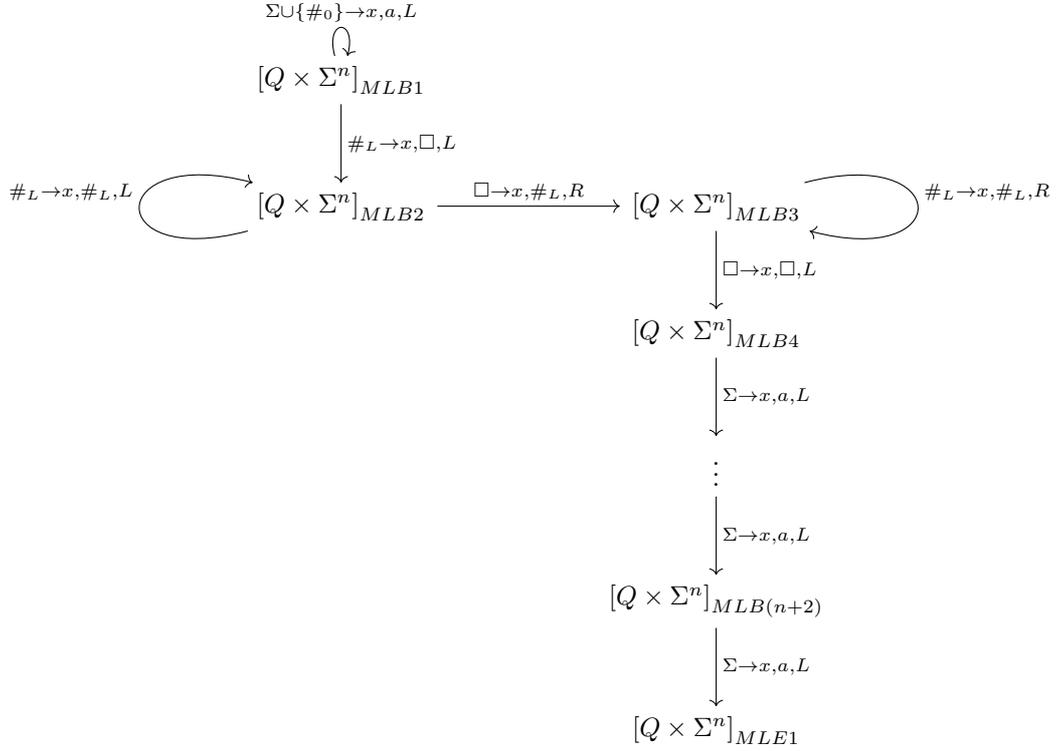
$$
\begin{tikzcd}[row sep = 3em, column sep = 6em]
{[Q \times \Sigma^n]}_{MLB1} \arrow[loop above, "{\Sigma \cup \{\#_0\} \rightarrow x,a,L}"] \arrow{d}{\#_L \rightarrow x,\square,L} \\
{[Q \times \Sigma^n]}_{MLB2} \arrow[loop left, "{\#_L \rightarrow x,\#_L,L}"] \arrow{r}{\square \rightarrow x,\#_L,R} & {[Q \times \Sigma^n]}_{MLB3} \arrow[loop right, "{\#_L \rightarrow x,\#_L,R}"] \arrow{d}{\square \rightarrow x,\square,L} \\
  & {[Q \times \Sigma^n]}_{MLB4} \arrow{d}{\Sigma \rightarrow x,a,L} \\
  & \vdots \arrow{d}{\Sigma \rightarrow x,a,L} \\
  & {[Q \times \Sigma^n]}_{MLB(n+2)} \arrow{d}{\Sigma \rightarrow x,a,L} \\
  & {[Q \times \Sigma^n]}_{MLE1}
\end{tikzcd}
$$

Here, the objective is to shift the left border of the bottom row one column to the left, to make room for the outward flow of non-empty squares. This outward flow occurs at a rate no faster than one column per simulation step, so it is sufficient to displace the border one column outwards in this way each cycle. In a more typical construction, the borders will be shifted depending on whether the simulated tape head moves to the edge of the tape encoding, however such an approach may introduce ambiguity in the tape head movement of $SIM$, so we avoid it. 

After executing the shift, the tape head of $SIM$ returns to the square corresponding to $\bs{y}_L^{(n)}$. \\

\emph{Left edge case} \\

In the following, we refer to the function $\Psi$ as defined at the end of section 2.1.

$$
\begin{tikzcd}[row sep = 3em, column sep = 8em]
{[Q \times \Sigma^n]}_{MLE1} \arrow{d}{\Sigma \rightarrow (x,a),a,R} \\
{[Q \times \Sigma^n]}_{MLE2} \arrow{r}{\Sigma \rightarrow x,a,R} & \cdots \arrow{r}{\Sigma \rightarrow x,a,R} & {[Q \times \Sigma^n \times \Sigma]}_{MLE(n+1)} \arrow{d}{\Sigma \rightarrow (x,a),a,L} \\
{[Q \times \Sigma^n \times \Sigma^2]}_{MLE(2n+1)} \arrow{d}{\Sigma \rightarrow x,\Psi^{(n)}(x_1,\cdots,x_{n+1},\square,x_{n+2},x_{n+3}),R} & \cdots \arrow{l}{\Sigma \rightarrow x,a,L} & {[Q \times \Sigma^n \times \Sigma^2]}_{MLE(n+2)} \arrow{l}{\Sigma \rightarrow x,a,L} \\
{[Q \times \Sigma^n \times \Sigma^2]}_{MML1}
\end{tikzcd}
$$

Here, $SIM$ computes the superposition at the left most encoded square of the current row. Since the square to the left of this on the corresponding tape of $MULTI$ is not explicitly encoded, but assumed to be blank, this is an edge case. At section $MLE1$, $SIM$ loads $\bs{y}_L^{(n)}$ into its local state distribution, which transforms as $\bs{q} \otimes \bs{y}_0 \mapsto \bs{q} \otimes \bs{y}_0 \otimes \bs{y}_L^{(n)}$. The tape head then moves right until it is over $\bs{y}_{L+1}^{(n)}$, which is also loaded, producing the local state distribution $\bs{q} \otimes \bs{y}_0 \otimes \bs{y}_L^{(n)} \otimes \bs{y}_{L+1}^{(n)}$. The tape head returns to $\bs{y}_L^{(n)}$, and writes the superposition given by $\Delta \Psi (\bs{q} \otimes \bs{y}_0 \otimes \square \otimes \bs{y}_L^{(n)} \otimes \bs{y}_{L+1}^{(n)})$. \\

\emph{Main loop} \\

We again refer to the function $\Psi$. Note that the sections omitted in the ellipsis will also be equipped with the same loop transitions of their source and target sections. Most of these loops will be redundant, depending on the current row being computed.

$$
\begin{tikzcd}[row sep = 3em, column sep = 8em]
& & {[Q \times \Sigma^n \times \Sigma^2]}_{MRB1} \\
{[Q \times \Sigma^n \times \Sigma^2]}_{MML1} \arrow[loop above, "{\#_0 \rightarrow x,a,R}"] \arrow{r}{\Sigma \rightarrow x,a,R} & \cdots \arrow{r}{\Sigma \rightarrow x,a,R} & {[Q \times \Sigma^n \times \Sigma^2]}_{MML(2n)} \arrow[loop above, "{\#_0 \rightarrow x,a,R}"] \arrow{d}{\Sigma \rightarrow (x,a),a,L} \arrow[u, swap, bend right = 80, "{\#_R \rightarrow x,\square,R}"] \\
{[Q \times \Sigma^n \times \Sigma^3]}_{MML(3n)} \arrow[loop below, "{\#_0 \rightarrow x,a,L}"] \arrow{u}[swap]{\Sigma \rightarrow (x_1,\cdots,x_{n+1},x_{n+3},x_{n+4}),\Psi^{(n)}(x),R} & \cdots \arrow{l}{\Sigma \rightarrow x,a,L} & {[Q \times \Sigma^n \times \Sigma^3]}_{MML(2n+1)} \arrow[loop below, "{\#_0 \rightarrow x,a,L}"] \arrow{l}{\Sigma \rightarrow x,a,L} \\
\end{tikzcd}
$$

At $MML1$, $SIM$ has local state distribution $\bs{q} \otimes \bs{y}_0 \otimes \bs{y}_{i-1}^{(n)} \otimes \bs{y}_{i}^{(n)}$, for some $L < i < R$. From $MML1$ to $MML(2n)$, the tape head moves rightward from $\bs{y}_{i}^{(n)}$ until $\bs{y}_{i_1}^{(n)}$ is reached, whereupon it is loaded into the local state distribution. From $MML(2n+1)$ to $MML(3n)$, the tape head returns to $\bs{y}_{i}^{(n)}$, where it writes the superposition as given by $\Delta \Psi (\bs{q} \otimes \bs{y}_0 \otimes \bs{y}_{i-1}^{(n)} \otimes \bs{y}_{i}^{(n)} \otimes \bs{y}_{i+1}^{(n)})$. At the same time, the distribution $\bs{y}_{i-1}^{(n)}$ is no longer needed by $SIM$, so is dumped from the local state distribution, which transforms as $\bs{q} \otimes \bs{y}_0 \otimes \bs{y}_{i-1}^{(n)} \otimes \bs{y}_{i}^{(n)} \otimes \bs{y}_{i+1}^{(n)} \mapsto \bs{q} \otimes \bs{y}_0 \otimes \bs{y}_{i}^{(n)} \otimes \bs{y}_{i+1}^{(n)}$. The loop exits when $\#_R$ is encountered instead of a symbol from $\Sigma$, triggering the right border shift. \\

\emph{Right border shift and right edge case} \\

$$
\begin{tikzcd}[row sep = 3em, column sep = 8em]
{[Q \times \Sigma^n \times \Sigma^2]}_{MRB1} \arrow{r}{\square \rightarrow x,\square,R} & \cdots \arrow{r}{\square \rightarrow x,\square,R} & {[Q \times \Sigma^n \times \Sigma^2]}_{MRBn} \arrow{d}{\square \rightarrow x,\#_R,L} \\
{[Q \times \Sigma^n \times \Sigma^2]}_{MRB(2n-1)} \arrow{d}{\Sigma \cup \{\#_R\} \rightarrow x,\square,L} & \cdots \arrow{l}{\Sigma \cup \{\#_R\} \rightarrow x,\square,L} & {[Q \times \Sigma^n \times \Sigma^2]}_{MRB(n+1)} \arrow{l}{\Sigma \cup \{\#_R\} \rightarrow x,\square,L} \\
{[Q \times \Sigma^n \times \Sigma^2]}_{MRE1} \arrow{d}{\Sigma \rightarrow (x_1,\cdots,x_{n+1}),\Psi^{(n)}(x,\square),L} \\
{[Q \times \Sigma^n]}_{MRE2}\arrow[loop below, "{\Sigma \rightarrow x,a,L}"] \arrow{r}{\#_0 \rightarrow x,\#_0,L} & \cdots
\end{tikzcd}
$$

The analysis here mirrors the left edge case and border shift, so we omit further explanation. The tape head returns to the $\#_0$ in the bottom row and the remaining rows are computed. \\

\underline{State update phase} \\

We assume the tape has returned to the square formerly containing $\bs{y}_0^{(1)}$. All that remains is to update the simulated state distribution, which is achieved as follows.

$$
\begin{tikzcd}[column sep = 8em]
{[Q \times \Sigma^n]}_{S} \arrow{r}{\Sigma \rightarrow \delta_1(x),a,S} & {[Q]}_{R1}
\end{tikzcd}
$$

The local state distribution transforms as $\bs{q} \otimes \bs{y}_0 \mapsto \Delta \delta_1 (\bs{q} \otimes \bs{y}_0)$. $SIM$ returns to the first section of the read phase, and the cycle is complete.

\end{proof}

\subsection{Universal Turing machine}

In this section, we introduce a design for a smooth relaxation preserving \emph{pseudo}-UTM. That is, a machine which simulates only machines with a maximal state count and tape alphabet size.

Our UTM will be a 2-tape machine. The notation generalises in the natural way as follows, with $a$ and $b$ denoting the read symbols of each tape and the components of $f$ mirroring those of $\delta$.

$$
\begin{tikzcd}[column sep = 24em]
{[X_i]}_i \arrow{r}{\Sigma_0^{(1)} \times \Sigma_0^{(2)} \rightarrow f_1(x,a,b), f_2(x,a,b), f_3(x,a,b), f_4(x,a,b), f_5(x,a,b)} & {[X_j]}_j
\end{tikzcd}
$$

We will also encounter a situation in which the state distribution is spread over two sections. In this case, the context will be a disjoint union of two contexts, say, $X \sqcup Y$. In analysing the smooth relaxation, we invoke the isomorphism $\mathbb{R}X \oplus \mathbb{R}Y \xrightarrow{\sim} \mathbb{R}(X \sqcup Y)$ sending $x \oplus 0 \mapsto x$ and $0 \oplus y \mapsto y$, thus denote local state distributions as direct sums.

\begin{theorem}
There exists a smooth relaxation preserving psuedo-UTM.
\end{theorem}

\begin{proof}
Let $M$ be a single tape TM with states $Q$, alphabet $\Sigma$ and transition function $\sigma$. Let $U$ be a $2$-tape machine. The tape alphabet $\Sigma_U$ of $U$ shall contain $Q \sqcup \Sigma$ along with an auxiliary symbol $\#$. Let $\bs{y}_i$ be the symbol distribution at index $i$ on the tape of $M$. The first tape will contain a sequence of tuples bordered by $\#$, encoding $M$. Each tuple will be of the form $(q,\sigma,\delta_1(q,\sigma),\delta_2(q,\sigma),\delta_3(q,\sigma))$, for a choice of $q \in Q$ and $\sigma \in \Sigma$. The tape head will be positioned over the left $\#$. The second tape will be identical to the tape of $M$. We shall call the first tape the \emph{description tape} and the second tape the \emph{working tape}.

The full construction of $U$ is given as follows.

$$
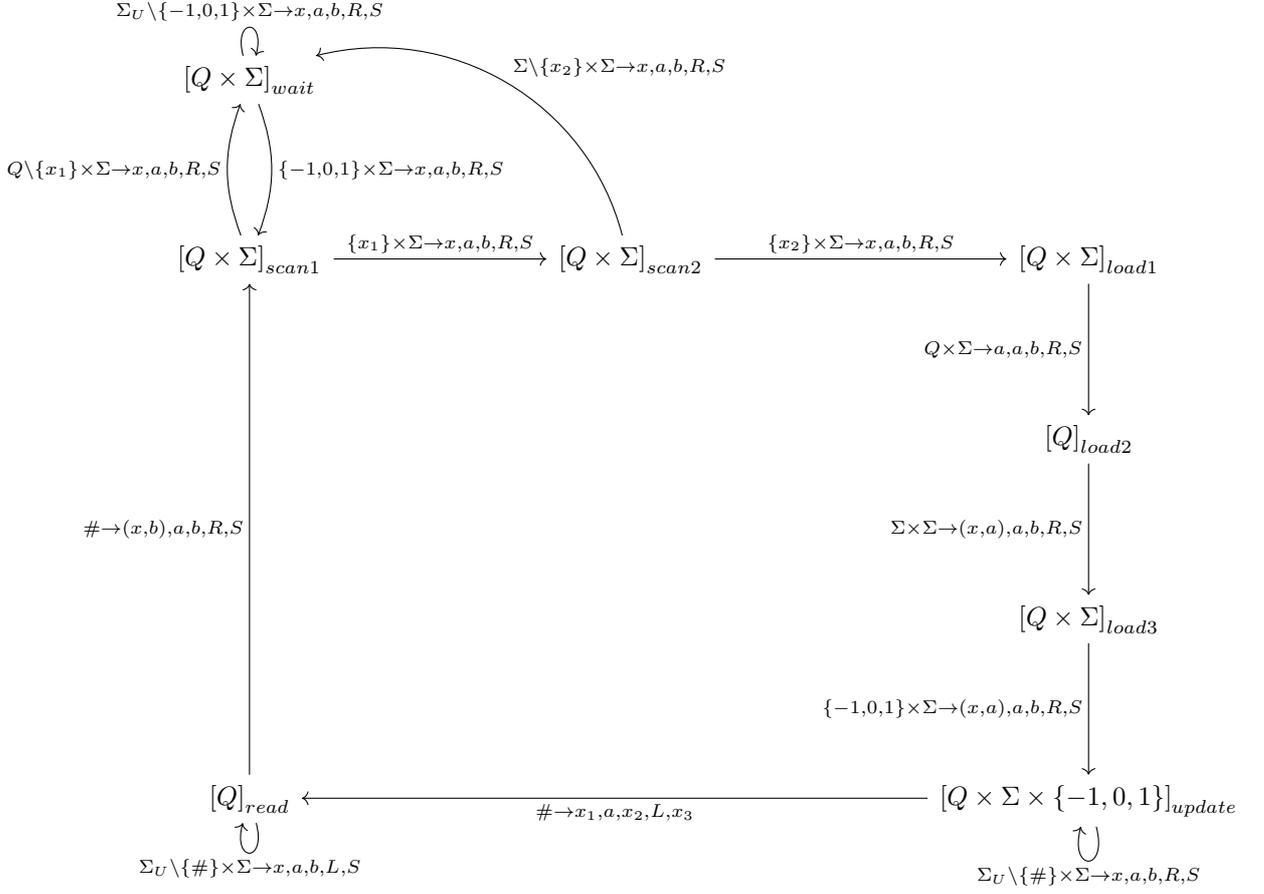
\begin{tikzcd}[row sep = 5em, column sep = 8em]
{[Q \times \Sigma]}_{wait} \arrow[loop above, "{\Sigma_U \setminus \{-1,0,1\} \times \Sigma \rightarrow x,a,b,R,S}"] \arrow[d, bend left = 20, "{\{-1,0,1\} \times \Sigma \rightarrow x,a,b,R,S}"] \\
{[Q \times \Sigma]}_{scan1} \arrow[u, bend left = 20, "{Q \setminus \{x_1\} \times \Sigma \rightarrow x,a,b,R,S}"] \arrow[r, "{\{x_1\} \times \Sigma \rightarrow x,a,b,R,S}"] & {[Q \times \Sigma]}_{scan2} \arrow[ul, swap, bend right = 45, "{\Sigma \setminus \{x_2\} \times \Sigma \rightarrow x,a,b,R,S}"] \arrow[r, "{\{x_2\} \times \Sigma \rightarrow x,a,b,R,S}"] & {[Q \times \Sigma]}_{load1} \arrow[d, swap, "{Q \times \Sigma \rightarrow a,a,b,R,S}"] \\
& & {[Q]}_{load2} \arrow[d, swap, "{\Sigma \times \Sigma \rightarrow (x,a),a,b,R,S}"] \\
& & {[Q \times \Sigma]}_{load3} \arrow[d, swap, "{\{-1,0,1\} \times \Sigma \rightarrow (x,a),a,b,R,S}"] \\
{[Q]}_{read} \arrow[loop below, "{\Sigma_U \setminus \{\#\} \times \Sigma \rightarrow x,a,b,L,S}"] \arrow[uuu,"{\# \rightarrow (x,b),a,b,R,S}"] & & {[Q \times \Sigma \times \{-1,0,1\}]}_{update} \arrow[loop below, "{\Sigma_U \setminus \{\#\} \times \Sigma \rightarrow x,a,b,R,S}"] \arrow[ll,"{\# \rightarrow x_1,a,x_2,L,x_3}"] \\
\end{tikzcd}
$$

At the beginning of a simulation cycle, $U$ has a state distributed over section $read$, with local distribution $\bs{q} \in \Delta Q$ mirroring the state distribution of $M$. Via the upwards tract, the read symbol distribution is loaded, yielding a local distribution $\bs{q} \otimes \bs{y}_0$. We may decompose this as
$$
\sum_{q,\sigma} \langle q,\bs{q} \rangle \langle \sigma,\bs{y}_0 \rangle q \otimes \sigma,
$$
where $q$ and $\sigma$ range over the standard bases of $\mathbb{R}Q$ and $\mathbb{R}\Sigma$. Henceforth, by \say{term}, we mean a summand in this decomposition. Meanwhile, the description tape head moves onto the beginning of the first tuple. From here, the objective of $U$ is to transport each term in the state distribution over section $scan1$ to an appropriate term in the target state distribution over section $update$. This will be clarified shortly. A term corresponding to local index $(q,\sigma)$ will arrive at local index $(\delta_1(q,\sigma),\delta_2(q,\sigma),\delta_3(q,\sigma)$. Note that The resulting local state distribution will not encode a conditionally independent joint distribution over $Q \times \Sigma \times \{-1,0,1\}$, however any dependence will be erased by the smooth relaxation in the leftward tract, during execution of which the working tape is updated.

We now verify that the transportation of the state distribution from section $scan1$ to $update$ behaves as desired. Let the first tuple correspond to the pair $(q_1,\sigma_1)$. As $U$ scans the initial state and read symbol on the tuple, the two rightward tracts have the effect of filtering out the correct term, with section $wait$ serving as a kind of sieve. The residue terms (left side of the final disjoint union below) remain in section $wait$ until $U$ has finished scanning the given tuple, upon which they return to $scan1$. When the first tuple is scanned, the local state distribution transforms as follows.
$$
\begin{tikzcd}
\sum_{q,\sigma} \langle q,\bs{q} \rangle \langle \sigma,\bs{y}_0 \rangle q \otimes \sigma \arrow[d, mapsto] \\
\sum_{q \neq q_1,\sigma} \langle q,\bs{q} \rangle \langle \sigma,\bs{y}_0 \rangle q \otimes \sigma \sqcup \sum_{q = q_1,\sigma} \langle q,\bs{q} \rangle \langle \sigma,\bs{y}_0 \rangle q \otimes \sigma \arrow[d, mapsto] \\
\sum_{(q,\sigma) \neq (q_1,\sigma_1)} \langle q,\bs{q} \rangle \langle \sigma,\bs{y}_0 \rangle q \otimes \sigma \sqcup \langle q_1,\bs{q} \rangle \langle \sigma_1,\bs{y}_0 \rangle q_1 \otimes \sigma_1
\end{tikzcd}
$$

The three downward tracts of $load1$, $load2$ and $load3$ distribute the right term over the target context, as $U$ loads the target state, write symbol and move directions. The right term transforms as follows (whilst the left term is stagnant).
\begin{align*}
\langle q_1,\bs{q} \rangle \langle \sigma_1,\bs{y}_0 \rangle q_1 \otimes \sigma_1 & \mapsto \langle q_1,\bs{q} \rangle \langle \sigma_1,\bs{y}_0 \rangle \delta_1(q_1,\sigma_1) \\
& \mapsto \langle q_1,\bs{q} \rangle \langle \sigma_1,\bs{y}_0 \rangle \delta_1(q_1,\sigma_1) \otimes \delta_2(q_1,\sigma_1) \\
& \mapsto \langle q_1,\bs{q} \rangle \langle \sigma_1,\bs{y}_0 \rangle \delta_1(q_1,\sigma_1) \otimes \delta_2(q_1,\sigma_1) \otimes \delta_3(q_1,\sigma_1) \\
\end{align*}

The remaining terms transform similarly. Once all terms have been transported, the tract bridging $update$ and $read$ executes the write and parallel move operation, and the tape head on the description tape returns to the left $\#$ before the cycle repeats.

\end{proof}

We briefly consider the behaviour of the above pseudo-UTM when the codes themselves are allowed to contain uncertainty. That is, suppose we allow the target state, write symbol and move direction in each tuple to take distributions in $\Delta Q$, $\Delta \Sigma$ and $\Delta \{-1,0,1\}$ respectively. This amounts to a generalisation of the transition function to a map, $\bs{\delta} : Q \times \Sigma \rightarrow \Delta Q \times \Delta \Sigma \times \Delta \{-1,0,1\}$. Extending the above analysis reveals that the step function for codes with uncertainty is obtained from Definition 2.2 by replacing all instances of the classical transition function with this generalised version. Indeed, this is a simplification of the behaviour exhibited by the staged pseudo UTM of \cite{gops}.

\section{Closing remarks}
As alluded to in the introduction, it is not clear whether similar results exist for simulating Turing machines of arbitrary alphabet size using alphabets of size $2$. According to condition 2 of Definition 2.4, one would need to devise a means of encoding distributions over finite sets of arbitrary size as sequences of Bernoulli distributions. One possibility is to employ a one-hot style encoding, with the probability of a given element in the sample space given by a single bit with uncertainty. However, a distribution encoded in this way cannot be loaded into the state distribution using the techniques demonstrated in section 3. If a certain kind of asymmetry is introduced into the encoding, it becomes possible to load the distribution, however (is seems) this is at the expense of being able to write the distributions encoded with the correct asymmetry back onto the tape.

Constructing smooth relaxation preserving simulations using a size $2$ alphabet is closely related to constructing smooth relaxation preserving (true) universal Turing machine, since the latter will inevitably involve variable length encodings of distributions over sets of arbitrary size.

\section{Acknowledgements}
The author is in gratitude to Dan Murfet for his patience, open-mindedness and invaluable feedback during the writing of this manuscript.

\bibliography{bib} 

\begin{thebibliography}{JCW21}

\bibitem[CM19]{dotm}
James Clift and Daniel Murfet.
\newblock Derivatives of turing machines in linear logic.
\newblock 2019.

\bibitem[JCW21]{gops}
Daniel~Murfet James~Clift and James Wallbridge.
\newblock Geometry of program synthesis.
\newblock 2021.

\bibitem[Rog96]{sutm}
Yurii Rogozhin.
\newblock Small universal turing machines.
\newblock 1996.

\end{thebibliography}
\bibliographystyle{alpha}

\end{document}